%% file: garstka_arxiv_clique_merging.tex
\pdfoutput=1
\documentclass[12pt]{article}

\usepackage{arxiv}
\input{options.tex}

\input{tabledata/result_variables.tex}
\usepackage{graphicx}
\title{A clique graph based merging strategy for decomposable SDPs}

\author{
  Michael Garstka\thanks{The authors are with the Department of Engineering Science, University of Oxford, UK}
\\
  \texttt{michael.garstka@eng.ox.ac.uk}
   \And
   Mark Cannon$^*$\\
  \texttt{mark.cannon@eng.ox.ac.uk}
   \And
 Paul Goulart$^*$\\
  \texttt{paul.goulart@eng.ox.ac.uk}
  }

\begin{document}

\maketitle
\input{sections/00_abstract.tex}

\input{sections/01_introduction.tex}

\input{sections/02_preliminaries.tex}

\input{sections/03_chordal_decomposition.tex}

\input{sections/04_clique_merging.tex}

\input{sections/05_benchmarks.tex}
\input{sections/06_conclusion.tex}

\section*{Acknowledgements}
We would like to thank Vidit Nanda and Heather \mbox{Harrington} for their helpful suggestions. M.\ Garstka is supported by the Clarendon Scholarship.

\bibliography{clique_bibliography_arxiv.bib}

\end{document}

%% file: options.tex
\usepackage{parskip}

\usepackage[small,labelfont={bf,sf}]{caption}

\usepackage{abstract}

\usepackage{hyperref} 
\hypersetup{colorlinks,
            citecolor=mdarkgreen,
            linkcolor=mmred}
\AtBeginDocument{}%
\AtBeginDocument{}%

\usepackage{amsmath}
\usepackage{amssymb}
\usepackage{amsthm}
\usepackage{mathtools}
\usepackage[ruled,vlined, linesnumbered]{algorithm2e}
\usepackage{siunitx}
\usepackage{xargs}
\usepackage[table]{xcolor}  
\usepackage{enumerate}
\usepackage{balance}

\usepackage{csquotes} 
\usepackage{csvsimple} 
\usepackage{booktabs} 
\usepackage[para,online,flushleft]{threeparttable}
\usepackage{tikz}
\usepackage{pgfplots}
\usepackage{pgfplotstable}
\usepackage{gincltex}
\usepackage{pgfmath,pgffor}

\pgfplotsset{compat=newest} 
\pgfplotsset{
      table/search path={figs/data},
    }

\usetikzlibrary{calc,arrows,intersections, matrix, fit, quotes, shapes, quotes}

\tikzset{
    sbox/.style={
        draw,
        rectangle split,
        rectangle split parts=2,
        text centered,
    },
    sbox+/.style={
        sbox,
        rectangle split every empty part={},
        rectangle split empty part width=\widthof{#1},
        rectangle split empty part height=\heightof{#1},
        rectangle split empty part depth=\depthof{#1},
    },
}

\newcommand{\abs}[1]{\left|#1\right|}
\newcommand{\mc}[1]{\mathcal{#1}}

\newtheorem{theorem}{Theorem}

\theoremstyle{definition}
\newtheorem{definition}{Definition}[section]

\theoremstyle{remark}

\makeatletter
\newcommand{\removelatexerror}{\let\@latex@error\@gobble}
\makeatother

\definecolor{mmred}{HTML}{e0024c}
\definecolor{mred}{rgb}{0.863,0.129,0.302}
\definecolor{mgreen}{rgb}{0,0.549,0}   
\definecolor{mblue}{rgb}{0,0.392,0.871} 
\definecolor{mpurple}{rgb}{0.4,0.0,0.4} 
\definecolor{morange}{rgb}{1,0.4,0} 
\definecolor{mteal}{rgb}{0.216,0.784,0.671} 
\definecolor{mdarkgreen}{RGB}{11, 91, 35}
\definecolor{oxfordblue}{RGB}{2,46,95}

\newcommand{\myCellColor}{red!25}

\newcommand{\OptProblem}[5][]{
  \begin{aligned}\label{#1}
    \begin{array}{ll}
    \underset{#3}{\rm{#2}}&\hspace{-0ex}#4\vspace{.5ex}\\
    \rm{subject~to}&#5
  \end{array}
\end{aligned}
}

\newcommand{\OptProblemExtra}[5][]{
  \begin{aligned}\label{#1}
    \begin{array}{lll}
    \underset{#3}{\rm{#2}}&\hspace{-0ex}#4\vspace{.5ex}\\
    \rm{subject~to}&#5
  \end{array}
\end{aligned}
}

\newcommand{\pkg}[1]{\texttt{#1}}

\renewcommand{\Re}{\mathbb{R}}
\newcommand{\Sym}[1][n]{\mathbb{S}^{#1}}
\newcommand{\Psd}[1][n]{\mathbb{S}_{+}^{#1}}
\newcommand{\clique}[1]{\mathcal{C}_{#1}}

%% file: tabledata/result_variables.tex
\newcommand{\gmeanRatios}{0.613}
\newcommand{\minRatio}{0.458} 
\newcommand{\minRatioName}{mcp500-3}

%% file: sections/00_abstract.tex
\begin{abstract}
Chordal decomposition techniques are used to reduce large structured positive semidefinite matrix constraints in semidefinite programs (SDPs).
The resulting equivalent problem contains multiple smaller constraints on the nonzero blocks (or cliques) of the original problem matrices.
This usually leads to a significant reduction in the overall solve time.
A further reduction is possible by remerging cliques with significant overlap.
The degree of overlap for which this is effective is dependent on the particular solution algorithm and hardware to be employed.
We propose a novel clique merging approach that utilizes the clique graph to identify suitable merge candidates.
We show its performance by comparing it with two existing methods on selected problems from a benchmark library.
Our approach is implemented in the latest version of the conic ADMM-solver \pkg{COSMO}.
\end{abstract}

%% file: sections/01_introduction.tex
We consider the primal-form semidefinite program (SDP):
\begin{equation}
\OptProblem[eq:primal_sdp]{minimize}{}{\langle C, X \rangle}{ \langle A_i, X \rangle = b_i, \quad i =1,\ldots, m\\ & X \in \Psd,}
\end{equation}
with variable $X$ and coefficient matrices $A_i, C \in \Sym$. The corresponding dual problem is
\begin{equation}
\OptProblem[eq:dual_sdp]{maximize}{}{b^\top y}{\displaystyle \sum_{i = 1}^{m}A_i y_i + S  = C\\ & S \in \Psd,}
\end{equation}
with dual variable $y \in \Re^m$ and slack variable $S$.
Semidefinite programming is used to solve problems that appear in a variety of applications such as portfolio optimisation, robust control, and optimal power flow problems.
Algorithms to solve SDPs, most notably interior point methods, have been known since the 1980s~\cite{Nesterov_1988}.
However, the recent trend to use models based on large quantities of data leads to SDPs whose dimensions challenge established solver algorithms.

Two main approaches are commonly used to deal with this challenge.
The first approach is to use first-order methods (FOMs) as in~\cite{ODonoghue_2016} or in~\cite{Zheng_2019}.
FOMs typically trade-off moderate accuracy solutions for a lower per-iteration computational cost and can therefore handle large problems more easily.

The second approach is to exploit sparsity in the problem data.
The authors in \cite{Fukuda_2001} showed that if the coefficient matrices $A_i, C$ exhibit an \textit{aggregate sparsity structure} represented by a \textit{chordal graph} $G(V, E)$, then the original primal and dual forms in~\eqref{eq:primal_sdp} and \eqref{eq:dual_sdp} can be decomposed.
These equivalent problems involve only positive semidefinite constraints on the nonzero blocks of the sparsity pattern which can lead to a significant reduction in the dimension of each constraint, thereby reducing solve time.
The equivalent primal problem is given by
\begin{equation}
\OptProblemExtra[eq:primal_decomp_sdp]{minimize}{}{\langle C, X \rangle}{ \langle A_i, X \rangle = b_i, & i =1,\ldots, m\\ & X_\ell = T_\ell X T^\top_\ell, & \ell = 1,\ldots,p \\ & X_\ell \in \Psd[|\mc{C}_\ell|], & \ell = 1,\ldots, p,}
\end{equation}
where the blocks $X_\ell$ are represented by subgraphs, called cliques, denoted $\mc{C}_\ell$.
Additional constraints using entry-selector matrices $T_\ell$, see~\eqref{eq:entry-selector}, enforce equality of the overlapping entries in $X$.
Following~\cite{Fukuda_2001} we refer to this conversion as the \textit{domain-space decomposition}.
The dual of this problem can be obtained by applying the \textit{range-space decomposition}:
\begin{equation}
\OptProblem[eq:dual_decomp_sdp]{maximize}{}{ b^\top y}{
\displaystyle \sum_{i = 1}^{m}A_i y_i + \sum_{\ell = 1}^p  T^\top_\ell S_\ell T_\ell = C\\
  & S_\ell \in   \Psd[|\mc{C}_\ell|],  \quad \ell = 1,\ldots, p.}
\end{equation}
Notice that the number and dimension of the block variables $X_\ell$ and $S_\ell$ depend only on the choice of cliques in the graph.
Starting from an initial decomposition we can merge two cliques $\mc{C}_i$ and $\mc{C}_j$ into a single clique with dimension $\abs{\mc{C}_i \cup \mc{C}_j}$.
Consequently, merging blocks has two opposing effects.
It increases the size of the blocks while decreasing the number of equality constraints.
Therefore, to evaluate the effect of a merge on the per-iteration time of a solver algorithm one has to take into account both the overlap between the cliques and the main linear algebra operations involved in each iteration.
\subsection*{Related work}\vspace{-1ex}
Heuristic methods to merge cliques have been proposed for interior-point methods.
The authors in~\cite{Nakata_2003} suggest traversing the \textit{clique tree}, a subset of all clique pairs with overlapping entries.
For each edge in the tree they merge corresponding cliques if the number of common entries relative to the cardinality of the individual cliques is higher than some threshold value, chosen heuristically to balance the block sizes and the number of additional equality constraints.
The methods are implemented in the \pkg{SparseCoLO} package~\cite{Fujisawa_2009}.

Similarly, the authors in~\cite{Sun_2014} suggest to traverse the clique tree and merge cliques if the amount of fill-in and the cardinality of the supernodes are below certain thresholds.
This approach is implemented in the \pkg{CHOMPACK} package~\cite{Andersen_2015}.

The authors in~\cite{Molzahn_2013} discuss clique merging in the context of a solver designed for large optimal power flow problems.
For each pair of adjacent cliques in the clique tree they determine how a merge would affect the problem dimension, i.e.\ the change in total number of variables and linking constraints.
They then greedily merge the blocks with the biggest reduction until the number of cliques decreases by a predefined percentage.

A limitation of existing methods is that they rely on heuristic parameters designed for a specific interior point implementation.
Furthermore, they consider only pairs of cliques that are adjacent in the clique tree.
We show in \autoref{sec:clique_merging} that two cliques with an advantageous merge are not necessarily adjacent in the clique tree.

With this paper we make the following contributions:
\begin{enumerate}[1.]
\item
We propose a clique merging strategy based on the clique intersection graph which considers all possible pair-wise merges and which can be tailored to the platform-specific matrix factorisation time.
\item
We use a weighting function for each pair of overlapping cliques that can be tailored to the specific algorithm used to solve the SDP. Specifically, we propose a weighting function for first-order methods that leverages the simpler relationship between clique sizes and per-iteration time, compared to interior-point methods. Consequently, we achieve consistently lower per-iteration times than with existing merging strategies.
\item
We provide a customisable implementation of our method in the latest version of the conic solver package \pkg{COSMO}~\cite{garstka_2019}.
\end{enumerate}

\subsection*{Outline}
In \autoref{sec:preliminaries} we define graph related concepts and describe how a graph-represented sparsity pattern can be used to decompose the primal and dual form of a SDP.
In \autoref{sec:clique_merging} we briefly outline two existing merging strategies based on the construction of the clique tree and describe our clique graph based approach.
In \autoref{sec:benchmarks} we then preprocess a number of benchmark problems using the different strategies and solve them with the same first order solver.
Consequently, we compare the impact of each strategy on the number of iterations, the per-iteration time, and the total solve time of the algorithm.
\autoref{sec:conclusion} concludes the paper.

%% file: sections/02_preliminaries.tex
\section{Graph Preliminaries}\label{sec:preliminaries}
In the following we define graph related concepts and how they relate to the sparsity structure of a matrix.
A good overview on this topic is provided by~\cite{Vandenberghe_2015}.
We consider the  \textit{undirected graph} $G(V, E)$ with vertex set $V$ and edge set $E \subseteq V \times V$.
Two vertices $v_1$, $v_2$ are \textit{adjacent} if $\{v_1, v_2 \} \in E$.
A \textit{cycle} is a path of edges (i.e.\ a sequence of distinct edges) joining a sequence of vertices in which only the first and last vertices are repeated.
A graph is called \textit{complete} if all vertices are pairwise adjacent.
We follow the convention of ~\cite{Vandenberghe_2015} by defining a \textit{clique} as a subset of vertices $\clique{} \subseteq V$ that induces a \textit{maximal} complete subgraph of $G$.

\begin{figure*}
\centering
  \input{figs/clique_graph_examples.tex}
  \caption{(a)~Aggregate sparsity pattern, (b)~sparsity graph $G(V,E)$, and (c)~clique tree $\mc{T}(\mc{B}, \mc{E})$.}
  \label{fig:example_structures}
\end{figure*}
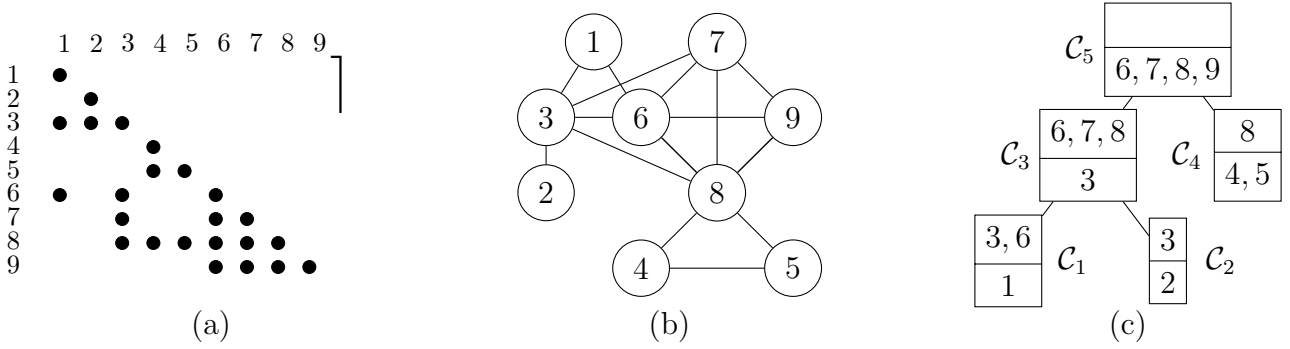

The decomposition theory described in Section~1b relies on a subset of graphs that exhibit the important property of \textit{chordality}.
A graph is \textit{chordal} (or \textit{triangulated}) if every cycle of length greater than three has a \textit{chord}, which is an edge between nonconsecutive vertices of the cycle.
A non-chordal graph can always be made chordal by adding extra edges.
An undirected graph with $n$ vertices can be used to represent the sparsity pattern of a symmetric matrix $S \in \Sym$.
Every nonzero entry $S_{ij} \neq 0$ in the lower (or upper) triangular part of the matrix introduces an edge $(i, j) \in E$.
An example of a sparsity pattern and the associated graph is shown in~\mbox{\hyperref[fig:example_structures]{Figure~\ref*{fig:example_structures}(a--b)}}.

For a given sparsity pattern $G(V, E)$, we define the following symmetric sparse matrix cones:
\begin{align*}
\Sym\left(E, 0\right) &\coloneqq \left\{ S \in \Sym \mid S_{ij} = S_{ji} = 0, \, \text{if}\, i \neq j, \, (i,j) \notin E\right\},\\
\Psd\left(E, 0\right) &\coloneqq \left\{ S \in \Sym(E, 0) \mid S \succeq 0\right\}.
\end{align*}
This means that for a matrix $S \in \Sym\left(E, 0\right)$ the diagonal entries $S_{ii}$ and the off-diagonal entries $S_{ij}$ with $(i,j)\in E$ may be zero or nonzero.
Moreover, we define the cone of positive semidefinite completable matrices:
\begin{equation*}
  \Psd(E, ?) \coloneqq \left\{ Y \mid \exists  \hat{Y} \in \Psd,  Y_{ij} = \hat{Y}_{ij}, \text{if}\, i=j \,\text{or}\, (i,j) \in  E \right\}.
\end{equation*}
For a matrix $Y \in  \Psd(E, ?)$ we can find a positive semidefinite completion by choosing appropriate values for all entries $(i, j)  \notin E$.
An algorithm to find this completion is described in~\cite{Vandenberghe_2015}.

An important structure conveying substantial information about the nonzero blocks of a matrix, or equivalently the cliques of a chordal graph, is the \textit{clique tree} (or \textit{junction tree}). For a chordal graph $G$ let $\mc{B} = \{\mc{C}_1, \ldots, \mc{C}_p \}$ be the set of cliques. The clique tree $\mc{T}(\mc{B},\mc{E})$ is formed by taking the cliques as vertices and by choosing edges from $\mc{E} \subseteq \mc{B} \times \mc{B}$ such that the tree satisfies the \textit{running-intersection property}:

\begin{definition}[Running intersection property]~\\
For each pair of cliques $\clique{i}$, $\clique{j} \in \mc{B}$, the intersection $\clique{i} \cap \clique{j}$ is contained in all the cliques on the path in the clique tree connecting $\clique{i}$ and $\clique{j}$.
\end{definition}
This property is also referred to as \textit{clique-intersection property} in~\cite{Nakata_2003} and \textit{induced subtree property} in~\cite{Vandenberghe_2015}.
For a given chordal graph, a clique tree can be computed using the algorithm described in~\cite{Pothen_1990}.

The clique tree for an example sparsity pattern is shown in Fig.~\ref{fig:example_structures}(c).
For a clique $\clique{\ell}$ we refer to the first clique encountered on the path to the root as its \textit{parent clique} $\clique{\mathrm{par}}$.
Conversely $\clique{\ell}$ is called the \textit{child} of  $\clique{\mathrm{par}}$.
If two cliques have the same parent clique we refer to them as \textit{siblings}.
For each clique define the functions $\mathrm{par}(\clique{\ell})$ and $\mathrm{ch}(\clique{\ell})$ that return its parent clique and its set of child cliques.
Note that each clique in Fig.~\ref{fig:example_structures}(c) has been partitioned into two sets.
The upper row represents the \textit{separators} $\eta_\ell = \clique{\ell} \cap \mathrm{par}(\clique{\ell})$, i.e.\ all clique elements that are also contained in the parent clique.
We call the sets of the remaining vertices shown in the lower rows the \textit{clique residuals} or \textit{supernodes} \mbox{$\nu_\ell = \clique{\ell} \setminus \eta_\ell$}. Keeping track of which vertices in a clique belong to the supernode and the separator is useful as the information is needed to perform a positive semidefinite completion.
For a set of vertices $V$, the \textit{power set} $\left\{ W \mid W \subseteq V \right\}$ is denoted as $2^{V}$.

%% file: figs/clique_graph_examples.tex
\begin{tabular}[b]{@{}cccc@{}}
\begin{minipage}{0.33\linewidth}
\includegraphics{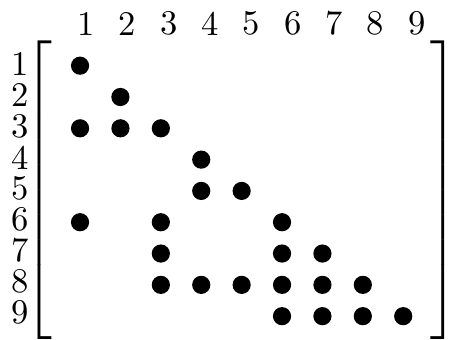}
\end{minipage}
&
\begin{minipage}{0.33\linewidth}
\centering
\includegraphics{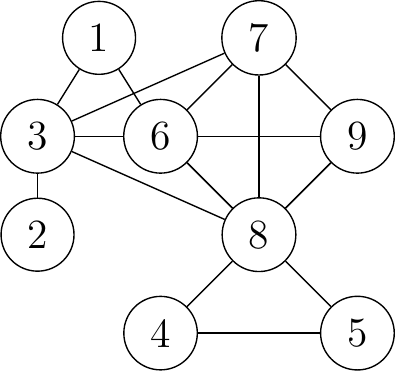}

 \end{minipage}
&
\begin{minipage}{0.33\linewidth}
\centering
\includegraphics{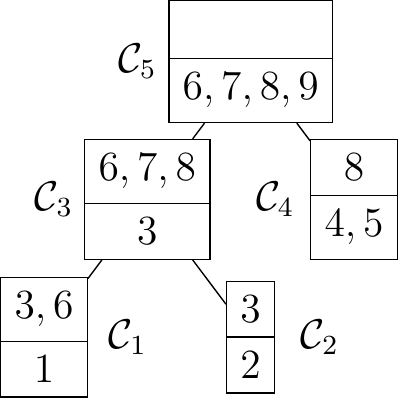}

\end{minipage}
\\
(a) & (b) & (c)
\end{tabular}

%% file: sections/03_chordal_decomposition.tex
\section*{1\MakeLowercase{b}. Chordal Decomposition}\label{sec:decomposition}
We next briefly describe how to apply chordal decomposition to an SDP. Let us assume that the problem matrices in~\eqref{eq:primal_sdp} and~\eqref{eq:dual_sdp} each have their own sparsity pattern
\begin{equation*}
  A_i \in \Sym[n](E_{A_i}, 0) \text{ and } C \in \Sym[n](E_{C}, 0).
\end{equation*}
The \textit{aggregate sparsity} of the problem is given by the graph $G(V, E)$ with edge set
\begin{equation*}
  E = E_{A_1} \cup E_{A_2} \cup \cdots \cup E_{A_m} \cup E_{C}.
\end{equation*}
In general $G(V, E)$ will not be chordal, but a \textit{chordal extension} can be found by adding edges to the graph.
We denote the extended graph as $G(V,\bar{E})$.
Finding the minimum number of edges to make the graph chordal is an NP-complete problem~\cite{Yannakakis_1981}. Consider a matrix $M$ of ones corresponding to the edge set $E$.
A commonly used heuristic method to find an extension is first to apply a reordering with approximate minimum fill-in~\cite{Amestoy_1996}.
Afterwards, a symbolic Cholesky factorisation is applied to the reordered matrix.
The Cholesky factor $L$ then defines a chordal extension with edge set $\bar{E}$.

Given sparsity information of the problem we can modify the matrix constraints in~\eqref{eq:primal_sdp} and~\eqref{eq:dual_sdp} to the respective sparse positive semidefinite matrix spaces:
\begin{equation}
\label{eq:sparse_constraints}
  X \in \Psd(\bar{E}, ?) \; \text{and} \;  S \in \Psd(\bar{E}, 0).
\end{equation}
We further define the entry-selector matrices \mbox{$T_\ell \in \Re^{ |\clique{\ell}|  \times n }$} for a clique $\clique{\ell}$:
\begin{equation}
\label{eq:entry-selector}
  (T_\ell)_{ij} \coloneqq \begin{cases}
    1, & \text{if } \clique{\ell}(i) = j\\
    0,  & \mathrm{otherwise,}
  \end{cases}
\end{equation}
where $\clique{\ell}(i)$ is the $i$th vertex of $\clique{\ell}$.
We can express the constraints in~\eqref{eq:sparse_constraints} in terms of multiple smaller coupled constraints using the theorems by~\cite{Grone_1984} and~\cite{Agler_1988}.

\begin{theorem}[Grone's theorem]\label{thm:grone}
Let $G(V,\bar{E})$ be a chordal graph with a set of maximal cliques $\{ \clique{1},\ldots,\clique{p}\}$. Then $X \in \Psd(\bar{E}, ?)$ if and only if
\begin{equation}
  X_\ell = T_\ell X T_\ell^\top \in \Psd[|\clique{\ell}|], \quad \forall \ell=1,\ldots, p.
\end{equation}
\end{theorem}
Applying this theorem to~\eqref{eq:primal_sdp} while restricting $X$ to the positive semidefinite completable matrix cone as in~\eqref{eq:sparse_constraints} yields the decomposed problem in~\eqref{eq:primal_decomp_sdp}. For the dual problem we utilise Agler's theorem, which is the dual to Thm.~\ref{thm:grone}:
\begin{theorem}[Agler's theorem]
Let $G(V,\bar{E})$ be a chordal graph with a set of maximal cliques $\{ \clique{1},\ldots,\clique{p}\}$.
Then $S \in \mathbb{S}^n_+(\bar{E},0)$ if and only if there exist matrices $S_\ell \in \Psd[|\clique{\ell}|]$ for $\ell=1,\ldots,p$ such that
\begin{equation}
    S = \sum_{\ell=1}^p T_\ell^\top S_\ell T_\ell. \label{eq:agler}
\end{equation}
\end{theorem}
With this theorem, we transform the dual form SDP in~\eqref{eq:dual_sdp} with the restriction on $S$ in~\eqref{eq:sparse_constraints} to arrive at~\eqref{eq:dual_decomp_sdp}.
Next, we show how to shape the sparsity pattern in the problem to reduce the per-iteration time of an SDP solver.

%% file: sections/04_clique_merging.tex
\section{Clique Merging}\label{sec:clique_merging}
Given an initial decomposition with edge set $\bar{E}$ and a set of cliques $\{\clique{1},\ldots,\clique{p}\}$, we are free to merge any number of cliques back into larger blocks.
This is equivalent to treating \textit{structural} zeros in the problem as \textit{numerical} zeros which leads to additional edges in the graph.
Looking at the decomposed problem in~\eqref{eq:primal_decomp_sdp} and~\eqref{eq:dual_decomp_sdp}, the effects of merging two cliques $\clique{i}$ and $\clique{j}$ are twofold:
\begin{enumerate}[1.]
  \item We replace two positive semidefinite matrix constraints of dimensions $\abs{\clique{i}}$ and $\abs{\clique{j}}$ with one constraint on a larger clique with dimension $\left| \clique{i} \cup \clique{j} \right|$, where the increase in dimension depends on the size of the overlap.
  \item We remove consistency constraints for the overlapping entries between $\clique{i}$ and $\clique{j}$, thus reducing the size of the linear system of equality constraints.
\end{enumerate}
When merging cliques these two factors have to be balanced.
The correct balance depends foremost on the used solver algorithm.
The authors in~\cite{Nakata_2003} and~\cite{Sun_2014} use the clique tree to search for favourable merge candidates. We will call these two approaches \texttt{SparseCoLO} and \textit{parent-child} strategy  in the following sections.
Before describing these methods, we define a procedure in \autoref{alg:merge_function} that describes how to merge a set of cliques within the set $\mc{B}$ and update the edge set $\mc{E}$ accordingly.
\begin{figure}[htb]
 \removelatexerror
\begin{algorithm}[H]
\SetKwInOut{Input}{Input}
\SetKwInOut{Output}{Output}
\SetKwRepeat{Do}{Do}{while}
\Input{A set of cliques $\mc{B}$ with edge set $\mc{E}$, a subset of cliques $\mc{B}_m = \{ \clique{m, 1} , \clique{m, 2},  \ldots, \clique{m, r}\}\subseteq \mc{B} $ to be merged.}
\Output{A reduced set of cliques $\hat{\mc{B}}$ with edge set $\hat{\mc{E}}$ and the merged clique $\clique{m}$.}
 $\hat{\mc{E}} \gets \mc{E} $\;
$\clique{m} \gets  \clique{m, 1} \cup   \clique{m, 2} \cup\cdots \cup \mc{C}_{m, r}$\;
$\hat{\mc{B}} \gets (\mc{B} \setminus \mc{B}_m) \cup \{\clique{m}\} $\;
  Remove edges  $ \{(\clique{i}, \clique{j}) \mid i \neq j, \; \clique{i}, \clique{j} \in \mc{B}_m \}$ in $\hat{\mc{E}}$\;
  Replace edges $\{(\clique{i}, \clique{j}) \mid \clique{i} \in \mc{B}_m,  \clique{j} \notin  \mc{B}_m \} $ with $(\clique{m}, \clique{j})$ in $\hat{\mc{E}}$\;
 \caption{Function $\mathrm{mergeCliques}(\mc{B}, \mc{E}, \mc{B}_m)$.}
   \label{alg:merge_function}
\end{algorithm}
\end{figure}
\subsection{Existing clique tree-based strategies}
The parent-child strategy described in~\cite{Sun_2014} traverses the clique tree in depth-first order and merges a clique $\clique{\ell}$ with its parent clique \mbox{$\clique{\mathrm{par}(\ell)} \coloneqq \mathrm{par}(\clique{\ell})$} if at least one of the two following conditions are met:
\begin{align}
   \left( \abs{\mc{C}_{\mathrm{par}(\ell)}} - \abs{\eta_\ell} \right) \left( \abs{\clique{\ell}} - \abs{\eta_\ell} \right) &\leq t_\mathrm{fill},\\
   \max \left\{ \abs{\nu_\ell}, \abs{\nu_{\mathrm{par}(\ell)}} \right\} &\leq t_\mathrm{size},
\end{align}
with heuristic parameters $t_{\mathrm{fill}}$ and $t_{\mathrm{size}}$.
The conditions keep the amount of extra fill-in and the supernode cardinalities below the specified thresholds.
The \texttt{SparseCoLO} strategy described in~\cite{Nakata_2003} and~\cite{Fujisawa_2006} considers parent-child as well as sibling relationships.
Given a parameter $\sigma > 0$, two cliques $\clique{i}, \clique{j}$ are merged if the following merge criterion holds
\begin{equation}
\label{eq:sparsecolo_criterion}
  \min \left\{ \frac{\abs{\clique{i} \cap \clique{j}}}{\abs{\clique{i}}}, \frac{\abs{\clique{i} \cap \clique{j}}}{\abs{\clique{j}}}  \right\} \geq \sigma.
\end{equation}
This approach traverses the clique tree depth-first, performing the following steps for each clique $\clique{\ell}$:
\begin{enumerate}[1.]
\item
For each clique pair \mbox{$\left\{ (\clique{i}, \clique{j}) \mid \clique{i}, \clique{j} \in \mathrm{ch}\left(\clique{\ell} \right) \right\}$}, check if~\eqref{eq:sparsecolo_criterion} holds, then:
\begin{itemize}
   \item $\clique{i}$ and $\clique{j}$ are merged, or
   \item if $(\clique{i} \cap \clique{j}) \supseteq \clique{\ell}$, then $\clique{i}$, $\clique{j}$, and $\clique{\ell}$ are merged.
 \end{itemize}
\item
For each clique pair $\left\{ \left(\clique{i}, \clique{\ell}\right) \mid \clique{i} \in \mathrm{ch}\left(\clique{\ell}\right) \right\}$, merge $\clique{i}$ and $\clique{\ell}$ if \eqref{eq:sparsecolo_criterion} is satisfied.
\end{enumerate}
We remark that the implementation of \pkg{SparseCoLO} follows the algorithm outlined here, but also employs a few additional heuristics.

\begin{figure}[htb]
  \centering
  \input{figs/nephew_uncle_example.tex}

  \caption{Sparsity graph (a) that can lead to a clique tree (b) with an advantageous \textquote{nephew-uncle} merge between $\clique{1}$ and $\clique{3}$.}
  \label{fig:nephew_merge}
\end{figure}
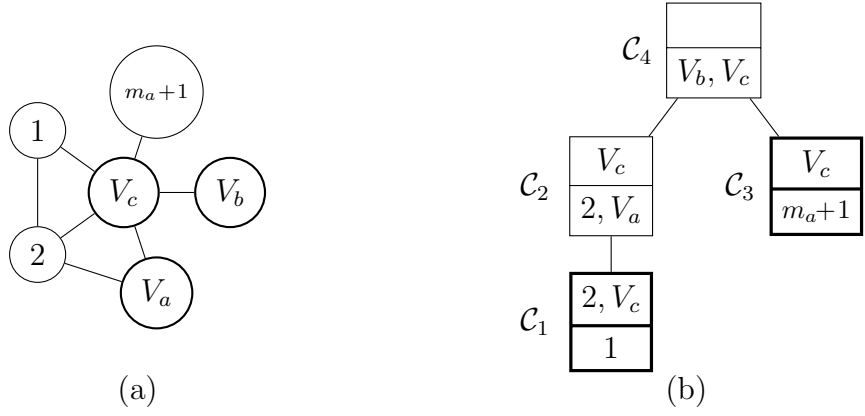

An advantage of the two approaches is that the clique tree can be computed easily and the conditions are inexpensive to evaluate.
However, a disadvantage is that choosing parameters that work well on a variety of problems and solver algorithms is difficult.
Secondly, in some cases it is beneficial to merge cliques that are not directly related on the clique tree.
To see this, consider a chordal graph $G(V,E)$ consisting of three connected subgraphs:
\begin{align*}
G_a(V_a, E_a),  \text{ with } V_a &= \{3, 4, \ldots, m_a\}, \\
G_b(V_b, E_b),  \text{ with } V_b &= \{m_{a}+2, m_{a}+3, \ldots, m_b\}, \\
G_c(V_c, E_c),  \text{ with } V_c &= \{m_{b}+1, m_{b}+2, \ldots, m_c\},
\end{align*}
and some additional vertices $\{1,2,m_a + 1\}$.
The graph is connected as shown in \hyperref[fig:nephew_merge]{Figure~\ref*{fig:nephew_merge}(a)}, where the complete subgraphs are represented as nodes $V_a,V_b,V_c$.
A corresponding clique tree is shown in \hyperref[fig:nephew_merge]{Figure~\ref*{fig:nephew_merge}(b)}.
By choosing the cardinality $\abs{V_c}$, the overlap between cliques $\clique{1} = \{1, 2\} \cup V_c$ and \mbox{$\clique{3} = \{m_a+1\} \cup V_c $} can be made arbitrarily large while $\abs{V_a}$, $\abs{V_b}$ can be chosen so that any other merge is disadvantageous.
However, neither the parent-child strategy nor \pkg{SparseCoLO} would consider merging $\clique{1}$ and $\clique{3}$ since they are in a \textquote{nephew-uncle} relationship.

\subsection{A new clique graph-based strategy}
To overcome the limitations of existing strategies we propose a merging strategy based on the \textit{clique-(intersection) graph} $\mc{G}(\mc{B}, \xi)$, where the edge set $\xi$ is defined as
\begin{equation*}
\xi = \left\{ (\clique{i}, \clique{j}) \mid  i \neq j, \; \clique{i}, \clique{j} \in \mc{B}, \; \abs{\clique{i} \cap \clique{j}} > 0 \right\} .
\end{equation*}
Let us further define an \textit{edge weighting function}\newline  \mbox{$e \colon 2^{V} \times 2^{V} \rightarrow \Re$}  that assigns a weight $w_{ij}$ to each edge  $(\clique{i}, \clique{j}) \in \xi$:
\[
e\left(\clique{i}, \clique{j}\right) = w_{ij}.
\]
This function is used to estimate the per-iteration computational savings of merging a pair of cliques depending on the targeted algorithm and hardware.
It is chosen to evaluate to a positive number if a merge would reduce the per-iteration time and to a negative number otherwise.
For a first-order method, whose per-iteration cost is dominated by an eigenvalue factorisation with complexity $\mc{O}\bigl(\abs{\clique{}}^3\bigr)$, a naive implementation would be:
\begin{equation}\label{eq:edge_weighting}
e(\clique{i}, \clique{j}) = \abs{\clique{i}}^3 + \abs{\clique{j}}^3 - \abs{\clique{i} \cup \clique{j}}^3.
\end{equation}
More sophisticated weighting functions can be determined empirically; see \autoref{sec:benchmarks}.
After a weight has been computed for each edge $(\clique{i}, \clique{j})$ in the clique graph, we merge cliques as outlined in \autoref{alg:cg_strategy}.
\begin{figure}[htb]
 \removelatexerror
\begin{algorithm}[H]
\SetKwInOut{Input}{Input}
\SetKwInOut{Output}{Output}
\SetKwRepeat{Do}{Do}{while}
\Input{A weighted clique graph $\mc{G}(\mc{B}, \xi)$. }
\Output{A merged clique graph $\mc{G}(\hat{\mc{B}}, \hat{\xi})$. }
$\hat{\mc{B}} \gets \mc{B}$ and $\hat{\xi} \gets \xi $\;
STOP $\gets$ false\;
\While{!STOP}{
  choose $\left(\clique{i}, \clique{j}\right)$ with maximum $w_{ij}$\;
  \eIf{$w_{ij} > 0$}{
    $\mc{B}_m \gets \{\clique{i}, \clique{j} \}$\;
    $\hat{\mc{B}}, \hat{\xi}, \clique{m} \gets \mathrm{mergeCliques}\left(\hat{\mc{B}}, \hat{\xi}, \mc{B}_m \right)$\;
    \For{each edge $(\clique{m}, \clique{\ell})   \in \hat{\xi}$}{
      update $w_{m \ell} \gets e(\clique{m}, \clique{\ell})$\;
    }
  }
  {
  STOP $\gets$ true\;
  }
}
 \caption{ Clique graph-based merging strategy.}
   \label{alg:cg_strategy}
\end{algorithm}
\end{figure}
Our strategy considers the edges in terms of their weights, starting with the clique pair $(\clique{i}, \clique{j})$ with the highest weight $w_{ij}$.
If the weight is positive, the two cliques are merged and the edge weights for all edges connected to the merged clique $\clique{m} = \clique{i} \cup \clique{j}$ are updated.
This process continues until no edges with positive weights remain.

The clique graph for the clique tree in \hyperref[fig:example_structures]{Figure~\ref*{fig:example_structures}(c)} is shown in \hyperref[fig:merged_example]{Figure~\ref*{fig:merged_example}(a)} with the edge weighting function in~\eqref{eq:edge_weighting}. Following \autoref{alg:cg_strategy} the edge with the largest weight is considered first and the corresponding cliques are merged, i.e.\ $\{3, 6, 7, 8\}$ and $\{6, 7, 8, 9\}$.
The revised clique graph $\mc{G}(\hat{\mc{B}}, \hat{\xi})$ is shown in \hyperref[fig:merged_example]{Figure~\ref*{fig:merged_example}(b)}.
Since no edges with positive weights remain, the algorithm stops.

\begin{figure}[htb]
  \centering
  \input{figs/merged_clique_graph.tex}
  \caption{(a)~Clique graph $\mc{G}(\mc{B}, \xi)$ of the clique tree in \hyperref[fig:example_structures]{Figure~\ref*{fig:example_structures}(c)}  with edge weighting function $e(\clique{i}, \clique{j}) = \abs{\clique{i}}^3 + \abs{\clique{j}}^3 - \abs{\clique{i} \cup \clique{j}}^3$ and (b)~clique graph  $\mc{G}(\hat{\mc{B}}, \hat{\xi})$ after merging the cliques $\{3, 6, 7, 8\}$ and $\{6, 7, 8, 9\}$ and updating edge weights. }
  \label{fig:merged_example}
\end{figure}

After \autoref{alg:cg_strategy} has terminated, it is possible to recompute a valid clique tree from the revised clique graph.
This can be done in two steps.
First, the edge weights in $\mc{G}(\hat{\mc{B}}, \hat{\xi})$ are replaced with new weights:
\begin{equation*}
  \tilde{w}_{ij} = \abs{\clique{i} \cap \clique{j}}, \; \text{for all } (\clique{i},\clique{j}) \in \hat{\xi}.
\end{equation*}
Second, a clique tree is then given by any \textit{maximum weight spanning tree} of the newly weighted clique graph, which can be computed using e.g.\ the algorithm described in~\cite{Kruskal_1956}.

Our merging strategy has some advantages over competing approaches.
Since the clique graph covers a wider range of merge candidates, it will consider edges that do not appear in clique tree-based approaches (such as the \textquote{nephew-uncle} example in \autoref{fig:nephew_merge}).
Moreover, the edge weighting function allows one to make a merge decision based on the particular solver algorithm and hardware used.
One downside is that this approach is more computationally involved than the other methods.
However, experiments show that the extra time spent on finding the clique graph, merging the cliques, and recomputing the clique tree is only a fraction of the total computational savings.

%% file: figs/nephew_uncle_example.tex
\begin{tabular}[b]{@{}cc@{}}
\begin{minipage}{.4\linewidth}
\tikzstyle{sbox}=[rectangle split,rectangle split parts=2,draw,text centered]
\centering
\includegraphics{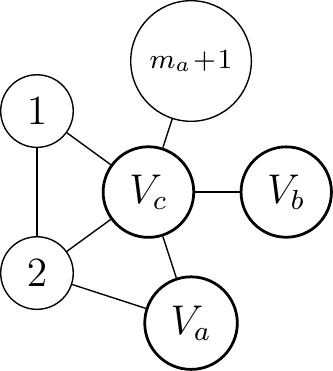}
\end{minipage}%
&
\begin{minipage}{.4\linewidth}
\centering
\includegraphics{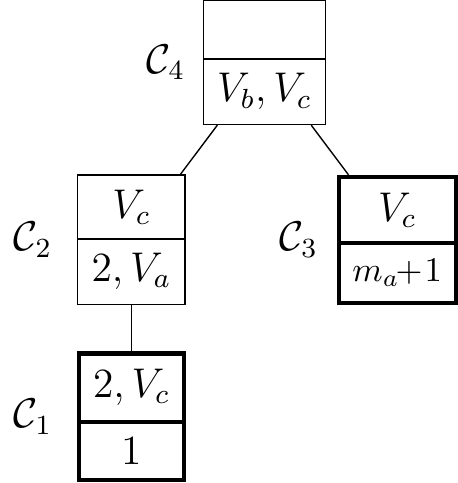}

\end{minipage}
\\
(a) & (b)
\end{tabular}

%% file: figs/merged_clique_graph.tex
\begin{tabular}[b]{@{}cc@{}}
\begin{minipage}{0.4\linewidth}
\centering
\includegraphics{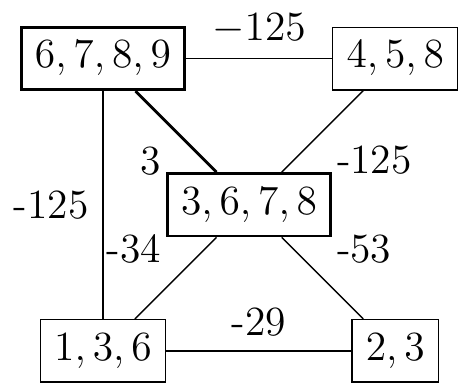}

\end{minipage}
&
\begin{minipage}{0.4\linewidth}
\centering
\includegraphics{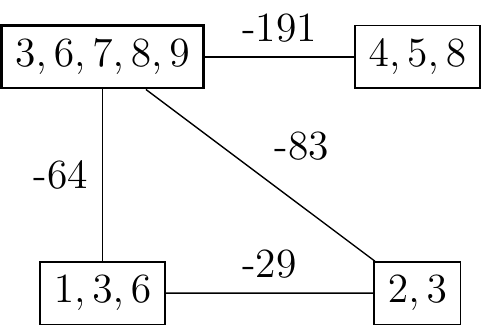}

\end{minipage}
\\
(a) & (b)
\end{tabular}

%% file: sections/05_benchmarks.tex
\section{Implementation and Results}\label{sec:benchmarks}
To compare the proposed merge approach with the clique tree-based strategies of~\cite{Nakata_2003} and~\cite{Sun_2014},
all three methods were used to preprocess sparse SDPs from SDPLib, a collection of SDP benchmark problems~\cite{Borchers_1999}.
Each strategy was given the same initial clique decomposition, and the resulting decomposed SDPs were solved using the first-order solver \pkg{COSMO}~\cite{garstka_2019}.
This section discusses how the different decompositions affect the per-iteration computation times of the solver.

\begin{figure}[htb]
\centering
\includegraphics{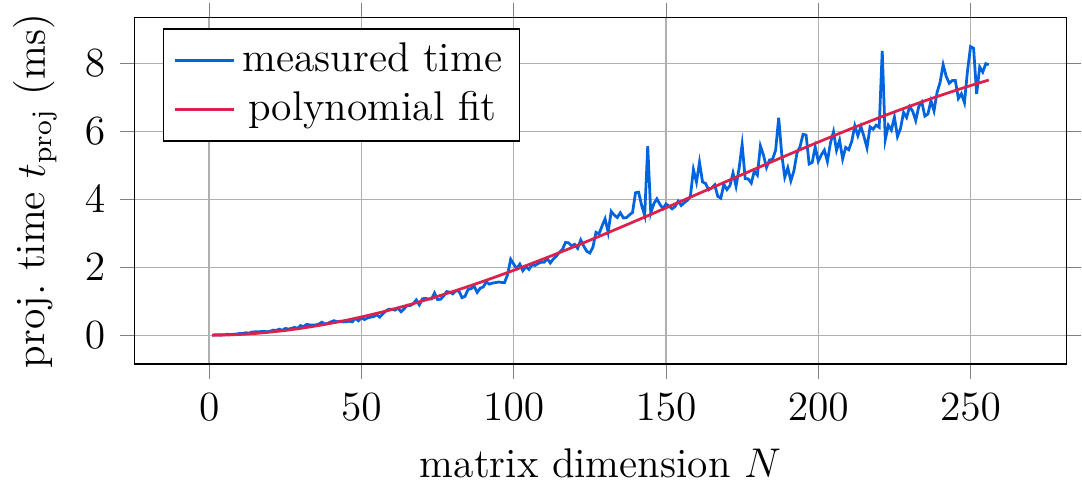}
\caption{Measured and estimated relationship between matrix size and execution time of the projection function in \pkg{COSMO}.}
  \label{fig:proj_estimate}
\end{figure}

For the strategy described in~\cite{Nakata_2003} we used the \pkg{SparseCoLO} package to decompose the problem.
The parent-child method by~\cite{Sun_2014} and our clique graph based method are available in the latest version of our conic solver \pkg{COSMO}.
For the former we chose the parameters $t_{\text{size}}=t_{\text{fill}}=9$.
We further investigate the effect of using different edge weighting functions.
Since \pkg{COSMO} is an ADMM-solver, the major operation affecting the per-iteration time is the projection step (see~\cite{garstka_2019} for more details).
This operation involves an eigenvalue decomposition of the matrices corresponding to the cliques.
Since the eigenvalue decomposition of a symmetric matrix of dimension $N$ has a complexity of $\mc{O}\left( N^3 \right)$, we define a \textit{nominal} edge weighting function as in~\eqref{eq:edge_weighting}.
However, the exact relationship will be different because the projection function involves copying of data and is affected by hardware properties such as cache size.
We therefore also consider an \textit{estimated} edge weighting function.
To determine the relationship between matrix size and projection time, the execution time of the relevant function inside \pkg{COSMO} was measured for different matrix sizes. We then approximated the relationship between projection time, $t_{\mathrm{proj}}$, and matrix size, $N$, as a polynomial:
\[
  t_{\mathrm{proj}}(N) = aN^3 + bN^2,
\]
where $a,b$ were estimated using least squares (\autoref{fig:proj_estimate}).
The estimated weighting function is then defined as
\begin{equation}
\label{eq:estimated_edge_weighting}
  e(\clique{i}, \clique{j}) = t_{\mathrm{proj}}\left( \abs{\clique{i}} \right) + t_{\mathrm{proj}} \left( \abs{\clique{j}} \right)  - t_{\mathrm{proj}} \left( \abs{\clique{i} \cup \clique{j}} \right).
\end{equation}

The merging strategies were compared for large, sparse SDP problems with chordal sparsity patterns from the SDPLib benchmark library. This problem set contains maximum cut problems, SDP relaxations of quadratic programs and Lovasz theta problems.
Six different cases were considered: no decomposition (\texttt{NoDe}), no clique merging (\texttt{NoMer}), decomposition using \pkg{SparseCoLO} (\texttt{SpCo}), parent-child merging (\texttt{ParCh}), and the clique graph-based method with nominal edge weighting (\texttt{CG1}) and estimated edge weighting (\texttt{CG2}).
All experiments were run on a MacBook with a \SI{2,6}{GHz} Intel Core i5-8259U CPU and \SI{8}{GB} of DDR3 RAM.
\pkg{COSMO} was configured to terminate with accuracy $\epsilon_{\mathrm{abs}}=\epsilon_{\mathrm{rel}}=\num{5e-4}$.
\pkg{SparseCoLO} was used with default parameters.
\autoref{tb:benchmark_results} shows the total solve time, the mean projection time, the number of iterations, the number of cliques after merging, and the maximum clique size of the sparsity pattern. The minimum value of each row is highlighted.
\begin{table}[htb]
\begin{center}
\caption{Benchmark results for different merging strategies.}\label{tb:benchmark_results}
\begin{threeparttable}
\resizebox{\textwidth}{!}{%
   \begin{tabular}{l l l l l l l l l l l l l}
     \toprule
     problem & \multicolumn{6}{c}{Solve time (\si{s}) } & \multicolumn{6}{c}{Projection time (\si{ms})}\\
         \cmidrule(lr){2-7}\cmidrule(lr){8-13}
      & NoDe\tnote{1} & NoMer\tnote{2} & SpCo\tnote{3} & ParCh\tnote{4} & CG1\tnote{5} & CG2\tnote{6} & NoDe & NoMer & SpCo & ParCh & CG1 & CG2 \\
     \midrule
     \input{./tabledata/benchmark_table_1.tex}
    \midrule
     problem & \multicolumn{6}{c}{Iterations} & \multicolumn{6}{c}{Number of cliques / Maximum clique size}\\
         \cmidrule(lr){2-7}\cmidrule(lr){8-13}
      & NoDe & NoMer & SpCo & ParCh & CG1 & CG2 &  NoDe & NoMer & SpCo & ParCh & CG1 & CG2 \\
     \midrule
    \input{./tabledata/benchmark_table_2.tex}
   \end{tabular}}
 \begin{tablenotes}
 \scriptsize{
    \item[1] no decomposition; \item[2] no merging; \item[3] \pkg{SparseCoLO} merging; \item[4] parent-child merging; \\
    \item[5] clique graph with nominal edge weighting~\eqref{eq:edge_weighting}; \item[6] clique graph with estimated edge weighting~\eqref{eq:estimated_edge_weighting}}
 \end{tablenotes}
 \end{threeparttable}
\end{center}
\end{table}

Our clique graph-based methods lead to a reduction in overall solve time. The method with estimated edge weighting function \texttt{CG2} achieves the lowest average projection times. The geometric mean of the ratios of projection time of \texttt{CG2} compared to the best non-graph method is $\gmeanRatios$, with a minimum ratio of $\minRatio$ for problem \mbox{\texttt{\minRatioName}}. Considering the number of cliques we see that \pkg{SparseCoLO} and \texttt{ParCh} merge more aggressively. Moreover, if the initial decomposition has a small maximum clique size, \pkg{SparseCoLO} seems to favor larger clique sizes. The merging strategies \texttt{ParCh}, \texttt{CG1} and \texttt{CG2} result in similar maximum clique sizes, with \texttt{CG1} being the most conservative in the number of merges.

%% file: tabledata/benchmark_table_1.tex
maxG11 & $29.7$ & $4.11$ & $7.9$ & $3.69$ & \cellcolor{\myCellColor} $\mathbf{2.72}$ & $2.82$ & $99.1$ & $15.3$ & $11.8$ & $12.3$ & $12.1$ & \cellcolor{\myCellColor} $\mathbf{9.2}$ \\ 
maxG32 & $320.98$ & $21.12$ & $27.08$ & $13.09$ & \cellcolor{\myCellColor} $\mathbf{12.47}$ & $15.79$ & $1105.0$ & $58.1$ & $57.8$ & $46.4$ & $38.3$ & \cellcolor{\myCellColor} $\mathbf{34.5}$ \\ 
maxG51 & $29.12$ & $28.04$ & $19.86$ & $9.59$ & \cellcolor{\myCellColor} $\mathbf{5.67}$ & $8.25$ & $171.4$ & $182.9$ & $191.9$ & $89.6$ & $54.3$ & \cellcolor{\myCellColor} $\mathbf{43.2}$ \\ 
mcp500-1 & $10.28$ & $1.04$ & $1.19$ & $0.78$ & $0.47$ & \cellcolor{\myCellColor} $\mathbf{0.37}$ & $40.4$ & $5.9$ & $6.9$ & $4.5$ & $3.4$ & \cellcolor{\myCellColor} $\mathbf{2.7}$ \\ 
mcp500-2 & $8.9$ & $10.25$ & $7.61$ & $5.97$ & $2.08$ & \cellcolor{\myCellColor} $\mathbf{1.95}$ & $35.2$ & $37.7$ & $34.5$ & $18.9$ & $11.6$ & \cellcolor{\myCellColor} $\mathbf{8.8}$ \\ 
mcp500-3 & $7.66$ & $22.69$ & $30.45$ & $15.76$ & $5.41$ & \cellcolor{\myCellColor} $\mathbf{4.35}$ & $35.7$ & $82.3$ & $83.9$ & $49.9$ & $24.2$ & \cellcolor{\myCellColor} $\mathbf{16.3}$ \\ 
mcp500-4 & $11.63$ & $51.37$ & $60.52$ & $21.92$ & \cellcolor{\myCellColor} $\mathbf{5.32}$ & $8.74$ & $39.6$ & $180.9$ & $132.3$ & $93.8$ & $36.8$ & \cellcolor{\myCellColor} $\mathbf{28.3}$ \\ 
qpG11 & $173.81$ & $6.05$ & $6.48$ & $7.65$ & $4.14$ & \cellcolor{\myCellColor} $\mathbf{3.87}$ & $397.2$ & $16.9$ & $11.9$ & $12.5$ & $13.4$ & \cellcolor{\myCellColor} $\mathbf{10.9}$ \\ 
qpG51 & $607.61$ & $138.38$ & $155.04$ & $150.14$ & $113.87$ & \cellcolor{\myCellColor} $\mathbf{85.19}$ & $749.9$ & $201.6$ & $185.6$ & $99.7$ & $58.8$ & \cellcolor{\myCellColor} $\mathbf{48.9}$ \\ 
thetaG11 & $225.89$ & $8.28$ & $37.16$ & $10.24$ & $9.01$ & \cellcolor{\myCellColor} $\mathbf{5.95}$ & $292.4$ & $20.9$ & $15.6$ & $14.8$ & $15.6$ & \cellcolor{\myCellColor} $\mathbf{12.8}$ \\ 
thetaG51 & $505.33$ & $82.48$ & $587.79$ & $103.47$ & \cellcolor{\myCellColor} $\mathbf{28.28}$ & $78.08$ & $193.4$ & $199.5$ & $204.6$ & $93.5$ & $58.8$ & \cellcolor{\myCellColor} $\mathbf{43.1}$ \\ 

%% file: tabledata/benchmark_table_2.tex
maxG11 & $280$ & $240$ & $640$ & $280$ & \cellcolor{\myCellColor} $\mathbf{200}$ & $280$ & $1/800$ & $598/24$ & $13/80$ & $207/32$ & $473/28$ & $411/38$ \\ 
maxG32 & $280$ & $320$ & $440$ & \cellcolor{\myCellColor} $\mathbf{240}$ & $280$ & $400$ & $1/2000$ & $1498/76$ & $21/210$ & $478/76$ & $1164/92$ & $468/126$ \\ 
maxG51 & $160$ & $120$ & $80$ & $80$ & $80$ & $160$ & $1/1000$ & $674/326$ & $181/322$ & $172/326$ & $448/362$ & $256/422$ \\ 
mcp500-1 & $240$ & $160$ & $160$ & $160$ & $120$ & $120$ & $1/500$ & $457/39$ & $451/44$ & $111/44$ & $437/54$ & $334/65$ \\ 
mcp500-2 & $240$ & $240$ & $200$ & $280$ & \cellcolor{\myCellColor} $\mathbf{160}$ & $200$ & $1/500$ & $363/138$ & $144/138$ & $111/140$ & $316/156$ & $223/177$ \\ 
mcp500-3 & $200$ & $240$ & $320$ & $280$ & $200$ & $240$ & $1/500$ & $259/242$ & $101/242$ & $70/242$ & $211/263$ & $134/301$ \\ 
mcp500-4 & $280$ & $240$ & $400$ & $200$ & \cellcolor{\myCellColor} $\mathbf{120}$ & $280$ & $1/500$ & $161/340$ & $63/346$ & $52/341$ & $105/368$ & $85/413$ \\ 
qpG11 & $400$ & $320$ & $520$ & $560$ & \cellcolor{\myCellColor} $\mathbf{280}$ & $320$ & $1/1600$ & $1398/24$ & $813/80$ & $296/32$ & $1273/28$ & $1211/38$ \\ 
qpG51 & $760$ & \cellcolor{\myCellColor} $\mathbf{600}$ & $720$ & $1360$ & $1800$ & $1640$ & $1/2000$ & $1674/326$ & $1182/304$ & $284/326$ & $1448/362$ & $1256/422$ \\ 
thetaG11 & $760$ & \cellcolor{\myCellColor} $\mathbf{360}$ & $2280$ & $640$ & $520$ & $400$ & $1/801$ & $598/25$ & $13/81$ & $207/33$ & $494/29$ & $423/41$ \\ 
thetaG51 & $2500$ & \cellcolor{\myCellColor} $\mathbf{320}$ & $2500$ & $920$ & $360$ & $1560$ & $1/1001$ & $676/324$ & $150/323$ & $169/324$ & $424/358$ & $202/425$ \\ 

%% file: sections/06_conclusion.tex
\section{Conclusion}\label{sec:conclusion}
A novel clique graph merging strategy to combine overlapping blocks in the aggregate sparsity pattern of structured SDPs is proposed.
The method considers all possible pair-wise merges and is customisable to the solver algorithm and hardware used.
An extension to our method would include information about the number of available CPU threads in the edge weighting function.
This would allow us to optimise the strategy for the parallel execution of the block-specific projection steps.
Benchmark tests show that our approach is able to reduce the projection time and the solve time of our first-order solver significantly compared to existing clique merging methods.